\documentclass[twoside]{article}
\title{Rings over which the
class of Gorenstein flat modules is closed under extensions}
\date{}
\usepackage[all]{xy}
\usepackage{amsfonts}
\usepackage{amsmath}
\usepackage{amssymb}
\usepackage{latexsym}
\usepackage{mathrsfs}
\usepackage{eufrak}

\newtheorem{thm}{\bf Theorem}[section]
\newtheorem{cor}[thm]{\bf Corollary}
\newtheorem{lem}[thm]{\bf Lemma}
\newtheorem{prop}[thm]{\bf Proposition}
\newtheorem{defn}[thm]{\bf Definition}
\newtheorem{defns}[thm]{\bf Definitions}
\newtheorem{rem}[thm]{\bf Remark}

\newtheorem{exmp}[thm]{\bf Example}

\catcode`\ç=13
\defç{\c{c}}
\catcode`\é=13
\defé{\'e}
\catcode`\à=13
\defà{\`a}
\catcode`\è=13
\defè{\`e}
\catcode`\â=13
\defâ{\^a}
\catcode`\ù=13
\defù{\`u}
\catcode`\ê=13
\defê{\^e}
\catcode`\î=13
\defî{\^\i}
\catcode`\ô=13
\defô{\^o}
\newcommand{\field}[1]{\mathbb{#1}}

\newcommand{\Q }{\field{Q}}
\newcommand{\Z }{\field{Z}}
\newcommand{\N }{\field{N}}

\def\proof{{\parindent0pt {\bf Proof.\ }}}

\def\fd{{\rm fd}}
\def\id{{\rm id}}

\def\Gfd{{\rm Gfd}}

\def\Im{{\rm Im}}

\def\Ker{{\rm Ker}}

\def\Ext{{\rm Ext}}
\def\Tor{{\rm Tor}}
\def\Hom{{\rm Hom}}

\def\wdim{{\rm wdim}}
\def\sup{{\rm sup}}

\newcommand{\GF}{\mathscr{GF}}

\newcommand{\cqfd}
{\hspace{1cm}
\rule{2mm}{2mm}%
\medbreak%
\par%
}

\begin{document}
\thispagestyle{empty}

\maketitle \vspace*{-1.5cm}

\begin{center}
{\large\bf Driss Bennis}

\bigskip

 Department of Mathematics, Faculty of Science and
Technology of Fez,\\ Box 2202, University S. M. Ben Abdellah Fez,
Morocco, \\[0.2cm]
 driss\_bennis@hotmail.com
\end{center}

\bigskip\bigskip
\noindent{\large\bf Abstract.} A ring $R$  is called left
GF-closed, if the class of all Gorenstein flat left $R$-modules is
closed under extensions. The class of left GF-closed rings
includes strictly the one of right coherent rings and the one of
rings of finite weak dimension.\\
In this paper, we investigate the Gorenstein flat dimension over
left GF-closed rings. Namely, we generalize the fact that the
class of all Gorenstein flat left modules is projectively
resolving over right coherent rings to left GF-closed rings. Also,
we generalize the characterization of Gorenstein flat left modules
(then of Gorenstein flat dimension of left modules) over right
coherent rings to left GF-closed rings. Finally,  using direct
products of rings, we show how to construct a left GF-closed ring
that is neither right coherent nor of finite weak
dimension.\bigskip

\small{\noindent{\bf Key Words.} Gorenstein flat dimension;
GF-closed rings; direct products of rings.}
\begin{section}{Introduction} Throughout the paper all rings
are associative with identity, and all modules are unitary.\\
Let $R$ be a ring and let $M$ be an $R$-module. The injective
(resp., flat) dimension of $M$ is denoted by $\id_R(M)$ (resp.,
$\fd_R(M)$). By $\wdim(R)$ we denote the weak dimension of $R$;
i.e., the supremum of the flat dimensions of all $R$-modules.\\
We say that $M$ is \textit{Gorenstein flat}, if there exists an
exact sequence of flat left $R$-modules, $$ \cdots\rightarrow
F_1\rightarrow F_0 \rightarrow F^0 \rightarrow F^1
\rightarrow\cdots,$$ such that  $M \cong \Im(F_0 \rightarrow F^0)$
and such that $I\otimes_R - $ leaves the sequence exact whenever
$I$ is an injective right $R$-module.\\
We use $\GF(R)$ to denote the class of all Gorenstein flat left
$R$-modules.\\
For a positive integer $n$, we say that $M$ has \textit{Gorenstein
flat dimension} at most $n$, and we write $\Gfd_R(M)\leq n$, if
$M$ has a Gorenstein flat resolution of length $n$; that is an
exact sequence:  $$ 0 \rightarrow G_n\rightarrow \cdots
\rightarrow G_0\rightarrow M \rightarrow 0,$$ where each $G_i$ is
Gorenstein flat left $R$-module (please see
\cite{LW,Rel-hom,HH}).\bigskip

The notion of Gorenstein flat modules was introduced and studied
over Gorenstein rings, by  Enochs, Jenda, and Torrecillas
\cite{GoPlat}, as a generalization of the notion of flat modules
in the sense that an $R$-module is flat if  and only if  it is
Gorenstein flat with finite flat dimension\,\footnote{\;  This
also holds over coherent rings (see \cite[note after Proposition
3.6]{BM}).}. In \cite{FCring}, Chen and Ding generalized known
characterizations of Gorenstein flat modules (then of the
Gorenstein flat dimension) over Gorenstein rings to $n$-FC rings
(coherent with finite self-FP-injective dimension). And in
\cite{HH}, Holm  relies on the use of character modules over
coherent rings to translate results for Gorenstein injective
modules to the setting of Gorenstein flat modules; and so he has
generalized the study of the Gorenstein flat dimension to coherent
rings.\bigskip

In this context, we enlarge the class of rings over which the
Gorenstein flat dimension is well-behaved. Namely, we investigate
the Gorenstein flat dimension over a new class of rings which we
call left GF-closed: a ring $R$ is called left GF-closed if the
class $\GF(R)$ of all Gorenstein flat left $R$-modules is closed
under extensions (Definition \ref{def-GF-closed}). The class of
left GF-closed rings includes strictly the one of right coherent
rings and also the one of rings of finite weak dimension (see the
note after Proposition \ref{pro-example-GF-closed}).\bigskip

The main results are in Section 2. First, we investigate the
behavior of Gorenstein flat modules in short exact sequences.
Recall the following:

\begin{defns}\label{def-resolv}
 \textnormal{Let $R$ be a ring and let $\mathscr{X}$ be a class of left
$R$-modules.\begin{enumerate}
    \item  We say that the class $\mathscr{X}$ is \textit{closed under
    extensions}, if for every short exact sequence of left $R$-modules $0\rightarrow A
    \rightarrow B \rightarrow C \rightarrow 0$, the condition $A$
and $C$ are in  $\mathscr{X}$ implies that $B$ is in
$\mathscr{X}$.
    \item We say that the class $\mathscr{X}$ is \textit{closed under
    kernels of epimorphisms}, if for every short exact sequence  of  left $R$-modules $0\rightarrow A
    \rightarrow B \rightarrow C \rightarrow 0$, the condition $B$
and $C$ are in  $\mathscr{X}$ implies that $A$ is in
$\mathscr{X}$.
    \item The class  $\mathscr{X}$  is said to be
     \textit{projectively resolving},
    if it contains all projective  left $R$-modules, and it is closed under
    both extensions and  kernels of epimorphisms; that is, for every
    short exact sequence  of  left $R$-modules $0\rightarrow A
    \rightarrow B \rightarrow C \rightarrow 0$ with  $C\in
    \mathscr{X}$, the conditions $A\in \mathscr{X}$ and  $B\in
    \mathscr{X}$ are equivalent.
\end{enumerate}}
\end{defns}

As a generalization of \cite[Theorem 3.7]{HH}, we show, over left
GF-closed rings, that the class $\GF(R)$ is  projectively
resolving and closed under direct summands (see Theorem
\ref{thm-close} and Corollary \ref{cor-GF-closed}). Then, we
generalize the characterization of Gorenstein flat dimension over
right coherent rings \cite[Theorem 3.14]{HH} to left GF-closed
rings. We end Section 2 with a generalization of \cite[Theorem
3.14]{HH}, which investigate the behavior of Gorenstein flat
dimension in short exact sequences (see Theorem
\ref{thm-short-Gfd})\bigskip

In Section 3, we study the Gorenstein flat dimension in direct products of
rings. The main result is the following (Theorem \ref{thm-Gflat-product}):\bigskip

\noindent\textbf{Theorem.}\textit{ Let
$R=\displaystyle\prod_{i=1}^n R_i$ be a direct product of rings
and let $N=N_1 \oplus \cdots\oplus N_n$  be a decomposition of a
left $R$-module $N$ into  left $R_i$-modules. Then,
$\Gfd_R(N)=\sup\{\Gfd_{R_i}(N_i)\, |\, i=0,...,n \}.$}\medskip

This enables us to construct a class of left GF-closed rings which
are neither right coherent nor of finite weak dimension (Example
\ref{exm-3.2}).

\end{section}\bigskip\bigskip


\begin{section}{Gorenstein flat dimension over left GF-closed rings}
In this section, we investigate the Gorenstein flat dimension
over the following kind of rings:

\begin{defn}\label{def-GF-closed} \textnormal{A ring $R$ is said to be
\textit{left GF-closed}, if $\GF(R)$ is closed under extensions.}
\end{defn}

Recall that a ring $R$ is called  \textit{right coherent}, if
every finitely generated right ideal $I$ is finitely presented;
that is there exists an exact sequence of right $R$-modules $ F_1
\rightarrow F_0 \rightarrow I \rightarrow 0$, where each $F_i$ is
finitely generated free.

\begin{prop}\label{pro-example-GF-closed} \begin{enumerate}
    \item  Every right coherent ring is left GF-closed.
    \item  Every ring with finite weak dimension is
left GF-closed.
\end{enumerate}
\end{prop}
\proof $1.$ Follows by \cite[Theorem 3.7]{HH}.\\
$2.$  Let $R$ be a ring such that $\wdim(R)\leq n$ for some
positive integer $n$. We show that the class of all Gorenstein
flat left $R$-modules and the class of all flat left $R$-modules
agree\,\footnote{\;  Compare to \cite[Corollary 3.8]{BM} in which
the rings are assumed to be commutative.}. Indeed, consider a
Gorenstein flat left $R$-module $G$. Then, by definition, there
exists an exact sequence of left $R$-modules $0\rightarrow
G\rightarrow F^0\rightarrow F^1\rightarrow \cdots \rightarrow
F^{n-1}\rightarrow N\rightarrow0$, where each $F^i$ is a flat left
$R$-module. By hypothesis, $\fd_R(N)\leq n$, and therefore $G$ is
flat, as desired.\cqfd\bigskip

It is well-known that there  exist right coherent (then left
GF-closed) rings with infinite weak dimension (see for example
\cite[Exercise 11 (f), p. 181]{Bou} or simply take a Gorenstein
ring, then Noetherian, with infinite weak dimension). Also, there
exist rings of finite weak dimension which are not right coherent.
This means that there exist left GF-closed rings which are not
right coherent (take for instance a non-semihereditary ring $R$
with $\wdim(R)\leq 1$, see for example \cite[Example 2.3]{KM}).
So, the class of left GF-closed rings includes strictly the one of
right coherent rings and the one of rings of finite weak
dimension. In the next section, we show that there exist a left
GF-closed ring which is neither coherent nor of finite weak
dimension (Example \ref{exm-3.2}). Namely, one would like to have
\textbf{every ring is left GF-closed}.\bigskip

We begin with the following fundamental result, which is, with
Corollary \ref{cor-GF-closed}, a generalization of   \cite[Theorem
3.7]{HH}.

\begin{thm}\label{thm-close} The following
conditions are equivalent for a ring $R$:
\begin{enumerate}
    \item $R$ is left GF-closed.
    \item  The class $\GF(R)$ is projectively resolving.
    \item  For every short exact sequence of left $R$-modules
$0\rightarrow G_1\rightarrow G_0\rightarrow M \rightarrow 0$,
where $ G_0$ and $G_1$ are Gorenstein flat. If $Tor^R_{1}(I,M)=0$
for all injective right $R$-modules $I$, then $M$ is Gorenstein
flat.
\end{enumerate}
\end{thm}

To prove this theorem, we need the following lemmas.

\begin{lem}\label{lem-close}  The following are equivalent for a left $R$-module $M$:
\begin{enumerate}
    \item  $M$ is Gorenstein flat.
    \item  $M$ satisfies the two following conditions: \begin{itemize}
        \item[(i)] $\Tor_i^R(I,M)=0$ for all $i>0$ and all
        injective right $R$-modules $I$, and
        \item[(ii)] There exists an exact sequence of left $R$-modules
        $0\rightarrow M \rightarrow
F^0 \rightarrow F^1 \rightarrow\cdots$, where each $F_i$ is flat,
such that  $I\otimes_R -$ leaves the sequence exact whenever $I$
is an injective right $R$-module.
    \end{itemize}
    \item There exists a short exact sequence of left
$R$-modules $ 0\rightarrow M\rightarrow F\rightarrow G \rightarrow
0,$ where $F$ is flat and $G$ is Gorenstein flat.
\end{enumerate}
\end{lem}
\proof Using the definition of Gorenstein flat modules, the
equivalence $(1) \Leftrightarrow (2)$ is  obtained by standard
argument. Also, by definition, we get immediately the implication
$(1)\Rightarrow (3)$.\\
We prove the implication $(3)\Rightarrow (2)$. Suppose that there
exists a short exact sequence of left $R$-modules: $$ (\alpha)=
\qquad0\rightarrow M\rightarrow F\rightarrow G \rightarrow 0,$$
where $F$ is flat and $G$ is Gorenstein flat. Let $I$ be an
injective right $R$-module. Since $G$ is Gorenstein flat and by
the equivalence $(1) \Leftrightarrow(2)$, $\Tor_{i+1}^R(I,G)=0$
for all $i\geq 0$. Then,  use the long exact sequence, $$
\Tor_{i+1}^R(I,G )  \rightarrow \Tor_{i}^R(I,M ) \rightarrow
 \Tor_{i}^R(I,F ),$$
to get $\Tor_{i}^R(I,M )=0$ for all $i>0$.\\
On the other hand, since $G$ is Gorenstein flat, there is an exact
sequence of left $R$-modules: $$ (\beta)= \qquad 0\rightarrow G
\rightarrow F^0 \rightarrow F^1 \rightarrow\cdots,$$  where each
$F_i$ is flat, such that  $I\otimes_R - $ leaves the sequence
exact whenever $I$ is an injective right $R$-module. Assembling
the sequences $(\alpha)$ and $(\beta)$, we get the following
commutative diagram:
$$\xymatrix @!0 @R=5mm  @C=1cm { 0\ar[r]& M \ar[r] & F
\ar[rr] \ar[rd] & & F^0  \ar[r]  & F^1  \ar[r]    & \cdots\\
 &  & & G \ar[rd]  \ar[ru]  &     &
    &  \\
     &  & 0  \ar[ru] &    & 0    &     &  }$$ such that
     $I\otimes_R  -  $ leaves the upper exact sequence exact whenever $I$ is an
injective right $R$-module, as desired.\cqfd

\begin{lem}\label{lem-close2}  Let $0\rightarrow
A\rightarrow B\rightarrow C \rightarrow 0$ be a short exact
sequence of left $R$-modules. If $A$ is Gorenstein flat and  $C$
is flat, then $B$ is Gorenstein flat.
\end{lem}
\proof Since $A$ is Gorenstein flat, there exists a short exact
sequence of left $R$-modules  $0\rightarrow A\rightarrow
F\rightarrow G \rightarrow 0$, where $F$ is flat and $G$ is
Gorenstein flat. Consider the following pushout diagram:
$$\xymatrix{
     &   0 \ar[d] & 0 \ar[d]  &   &  \\
0\ar[r]& A \ar[d] \ar[r] & B \ar@{-->}[d] \ar[r] & C \ar@{=}[d]
\ar[r] &
0\\
0\ar[r]& F \ar[d]  \ar@{-->}[r] & F'\ar[d] \ar[r] & C
\ar[r] & 0\\
 &  G \ar[d] \ar@{=}[r] &G \ar[d]  &   &  \\
     &   0 & 0   &   & }$$
In the sequence $0\rightarrow F\rightarrow F'\rightarrow C
\rightarrow 0$, both $F$ and $C$ are flat, hence so is $F'$. Then,
by the middle vertical sequence and from Lemma \ref{lem-close},
$B$ is Gorenstein flat, as desired.\cqfd\bigskip

In the following proof, $M^*$ stands for the character
module $\Hom_{\Z}(M, \Q/\Z)$ of a module $M$.\\

\noindent\textbf{Proof of Theorem \ref{thm-close}.}
$(1)\Rightarrow (2)$.  To claim that the class $\GF(R)$ is
projectively resolving, it suffices to prove that it is closed
under kernels of epimorphisms (see Definitions \ref{def-resolv}).
Then, consider a short exact sequence of left $R$-modules
$0\rightarrow A\rightarrow B\rightarrow C \rightarrow 0$, where
$B$  and $C$ are Gorenstein flat. We prove that $A$ is Gorenstein
flat. Since $B$ is Gorenstein flat, there exists a short exact
sequence of left $R$-modules $0\rightarrow B\rightarrow
F\rightarrow G\rightarrow 0$, where $F$ is flat and $G$ is
Gorenstein flat. Consider the following pushout diagram:
$$\xymatrix{
     &    &  0 \ar[d] & 0 \ar[d]  &  \\
0\ar[r]& A \ar@{=}[d] \ar[r] & B \ar[d] \ar[r] & C \ar@{-->}[d]
\ar[r] &
0\\
0\ar[r]& A   \ar[r] & F\ar[d] \ar@{-->}[r] & D \ar[d] \ar[r] &
0\\
 &   & G\ar[d] \ar@{=}[r] & G\ar[d]   &
 \\
 &   & 0 &0  & }$$ By the right vertical sequence and since $R$ is left GF-closed,
the $R$-module $D$ is Gorenstein flat. Therefore, by the middle
horizontal sequence and Lemma \ref{lem-close}, $A$ is Gorenstein
flat, as desired.\\ $(1)\Rightarrow (3)$. Since $G_1$ is
Gorenstein flat, there exists a short exact sequence of left
$R$-modules $0\rightarrow G_1\rightarrow F_1\rightarrow H
\rightarrow 0$, where $F_1$ is flat and $H$ is Gorenstein  flat.
Consider the following  pushout diagram:
$$\xymatrix{
     &   0 \ar[d] & 0 \ar[d]  &   &  \\
0\ar[r]& G_1 \ar[d] \ar[r] &G_0\ar@{-->}[d] \ar[r] & M \ar@{=}[d]
\ar[r] &
0\\
0\ar[r]& F_1 \ar[d]  \ar@{-->}[r] & D\ar[d] \ar[r] & M
\ar[r] & 0\\
 & H \ar[d] \ar@{=}[r] &H \ar[d]  &   &  \\
     &   0 & 0   &   & }$$
In the short exact sequence $ 0\rightarrow G_0\rightarrow
D\rightarrow H \rightarrow 0$ both $G_0$ and $H$ are Gorenstein
flat, then so is $D$ (since $R$ is left GF-closed). Then, there
exists a short exact sequence of left $R$-modules  $0\rightarrow
D\rightarrow F\rightarrow G \rightarrow 0$, where $F$ is flat and
$G$ is Gorenstein flat. Consider the following pushout diagram:
$$\xymatrix{
     &    &  0 \ar[d] & 0 \ar[d]  &  \\
0\ar[r]& F_1\ar@{=}[d] \ar[r] & D \ar[d] \ar[r] & M \ar@{-->}[d]
\ar[r] &
0\\
0\ar[r]&F_1 \ar[r] & F\ar[d] \ar@{-->}[r] & F'\ar[d] \ar[r] &
0\\
 &   & G\ar[d] \ar@{=}[r] & G\ar[d]   &
 \\
 &   & 0 &0  & }$$
We show that $F'$ is flat. Consider the sequence $0\rightarrow
M\rightarrow F'\rightarrow G\rightarrow 0$. Let $I$ be an
injective right $R$-module. By the exact sequence,
$$0=\Tor^R_1(I,M)\rightarrow \Tor^R_1(I,F')\rightarrow
\Tor^R_1(I,G)=0,$$ we get $$(*)\qquad\Tor^R_1(I,F')=0.$$ On the
other hand, consider the sequence $0\rightarrow F_1\rightarrow
F\rightarrow F'\rightarrow 0$. By \cite[Lemma 3.51]{Rot}, we have
the following short exact sequence of character modules:
$$ (\beta)=\qquad 0\rightarrow (F')^*\rightarrow F^*\rightarrow
(F_1)^* \rightarrow 0 .$$ From \cite[Theorem 3.52]{Rot}, $F^*$ and
$(F_1)^*$ are injective right $R$-modules. Then, by $(*)$ and from
\cite[Proposition 5.1, p. 120]{CE},
$$\Ext^1_R((F_1)^*, (F')^*)\cong (\Tor^R_1((F_1)^*,F'))^*=0.$$
Then, the sequence $ (\beta)$ splits, and so $(F')^* $ is
injective being a direct summand of the injective right $R$-module
$F^*$.
Therefore, $F'$ is a flat left $R$-module (by \cite[Theorem 3.52]{Rot}).\\
Finally, by Lemma \ref{lem-close} and the short exact sequence
$0\rightarrow M\rightarrow F'\rightarrow G\rightarrow 0$,
 $M$ is Gorenstein flat.\\
$(3)\Rightarrow (1)$.   Consider a short exact sequence of left
$R$-modules $0\rightarrow A\rightarrow B\rightarrow C \rightarrow
0$, where $A$  and $C$ are Gorenstein flat. We prove that $B$ is
Gorenstein flat. Let  $I$ be an injective right $R$-module.
Applying the functor $I\otimes_R -$ to the short exact sequence
$0\rightarrow A\rightarrow B\rightarrow C \rightarrow 0$, we get
the long exact sequence,
$$\Tor^R_i(I,A)\rightarrow \Tor^R_i(I,B )\rightarrow
\Tor^R_i(I,C ).$$ Then, $ \Tor^R_i(I,B)=0$ for all $i>0$ (since
$A$ and $C$ are Gorenstein flat and by Lemma \ref{lem-close}).\\
On the other hand, since $C$ is Gorenstein flat, there exists, by
definition, a short exact sequence of left $R$-modules
$0\rightarrow G\rightarrow F\rightarrow C \rightarrow 0$, where
$F$  is flat and $C$ is Gorenstein flat.  Consider the following
pullback diagram:
$$\xymatrix{
     &    &  0 \ar[d] & 0 \ar[d]  &  \\
 &   & G\ar[d] \ar@{=}[r] & G\ar[d]   &  \\
 0\ar[r]& A\ar@{=}[d] \ar[r] & D\ar@{-->}[d] \ar@{-->}[r] & F \ar[d]  \ar[r] & 0\\
0\ar[r]&A \ar[r] & B\ar[d] \ar[r]& C\ar[d] \ar[r] & 0\\
 &   & 0 &0  & }$$
Also, since $A$ is Gorenstein flat, there exists a short exact
sequence of  left $R$-modules $0\rightarrow A\rightarrow
F'\rightarrow  G' \rightarrow 0$, where $F'$ is flat and $G'$ is
Gorenstein  flat. Consider the following  pushout diagram:
$$\xymatrix{
     &   0 \ar[d] & 0 \ar[d]  &   &  \\
0\ar[r]& A\ar[d] \ar[r] & D\ar@{-->}[d] \ar[r] & F\ar@{=}[d]
\ar[r] &
0\\
0\ar[r]& F' \ar[d]  \ar@{-->}[r] & D'\ar[d] \ar[r] & F
\ar[r] & 0\\
 & G' \ar[d] \ar@{=}[r] &G' \ar[d]  &   &  \\
     &   0 & 0   &   & }$$
In the short exact sequence $ 0\rightarrow F'\rightarrow
D'\rightarrow F\rightarrow 0$ both $F'$ and $F$ are flat, then so
is $D'$. Then, by the short exact sequence $ 0\rightarrow
D\rightarrow D'\rightarrow G'\rightarrow 0$ and from Lemma
\ref{lem-close}, $D$ is Gorenstein flat. Finally, consider the
sort exact sequence $0\rightarrow G\rightarrow D\rightarrow B
\rightarrow 0$. We have $G$ and $D$ are Gorenstein flat, and
$\Tor^R_i( I,B)=0$ for all $i>0$ and all injective right
$R$-modules $I$. Therefore, by $(3)$, $B$ is Gorenstein flat. This
completes the proof.\cqfd

\begin{cor}\label{cor-GF-closed} If  $R$  is a left GF-closed ring, then
the class $\GF(R)$ is closed under direct summands.
\end{cor}
\proof  Use \cite[Propositions 1.4 and 3.2]{HH} and Theorem
\ref{thm-close}.\cqfd

\begin{rem}\label{rem-Gp} \textnormal{The Gorenstein
projective left modules have a similar characterization to the one
of the Gorenstein flat left modules in Lemma \ref{lem-close}.
Explicitly, a left $R$-module $M$ is Gorenstein projective if and
only if there exists a short exact sequence of left $R$-modules $
0\rightarrow M\rightarrow P\rightarrow G \rightarrow 0,$ where $P$
is projective and $G$ is Gorenstein projective. Thus, similar
argument to the one of the implication $(1) \Rightarrow (2)$ of
Theorem \ref{thm-close} gives a new proof to the fact that the
class of all Gorenstein projective left modules is projectively
resolving (see \cite[Theorem 2.5]{HH} and its proof).}
\end{rem}

Now, we give functorial descriptions of Gorenstein flat dimension
over left GF-closed rings. This  generalizes \cite[Theorem
3.14]{HH} which is proved for right coherent rings.

\begin{thm}\label{thm-GFD} Let $R$ be a ring and let  $M$ be a left $R$-module.
If $R$ is left GF-closed. Then, the following are equivalent for a
positive integer $n $:
\begin{enumerate}
    \item $\Gfd_R(M) \leq n$;
    \item $\Gfd_R(M) <\infty$ and $\Tor_{i}^R(I,M)=0$  for all $  i>n $  and all injective right $R$-modules
    $I$;
     \item  $\Gfd_R(M) <\infty$ and $\Tor_{i}^R(E,M)=0 $ for all  $i>n$ and all right $R$-modules
    $E$ with $\id_R(E)<\infty$;
    \item For every exact sequence of left $R$-modules
    $0\rightarrow K_{n}\rightarrow G_{n-1}\rightarrow \cdots \rightarrow G_0\rightarrow
    M\rightarrow 0$, if each $G_i $ is Gorenstein flat, then so is $K_{n}$.
\end{enumerate}
Furthermore, if $\Gfd_R(M) <\infty$, then
\begin{eqnarray*}
  \Gfd_R(M) &=& \sup\{ i\in \N \,|\,  \Tor_{i}^R(E,M)\not = 0 \ for \ some\ E\ with\ \id_R(E)< \infty\} \\
  &=& \sup\{ i\in \N \,|\, \Tor_{i}^R(I,M)\not = 0   \ for \ some\ injective\ right\ R\!-\!module\ I
  \}.
\end{eqnarray*}
\end{thm}

To prove this theorem, we need the following lemma.

\begin{lem}\label{lem-GFD} Let $M$ be a left $R$-module and consider
two exact sequences of left $R$-modules,
 $$\begin{array}{l}
 0 \rightarrow G_n\rightarrow G_{n-1}\rightarrow  \cdots \rightarrow G_0\rightarrow M \rightarrow
0,\ and\\
0 \rightarrow H_n\rightarrow H_{n-1}\rightarrow  \cdots
\rightarrow H_0\rightarrow M \rightarrow 0.
\end{array}$$
where $G_0,...,G_{n-1}$ and $H_0,...,H_{n-1}$ are Gorenstein flat.
If $R$ is left GF-closed, then $G_n$ is Gorenstein flat if  and
only if $H_n$ is Gorenstein flat.
\end{lem}
\proof Using Theorem \ref{thm-close}, Corollary
\ref{cor-GF-closed}, and \cite[Proposition 3.2]{HH}, the proof is
similar to the one of \cite[Theorem 1.2.7 $(i)\Rightarrow
(iii)$]{LW}.\cqfd\bigskip

\noindent\textbf{Proof of Theorem \ref{thm-GFD}.} First, note that
the last equalities follow immediately from the equivalences
between $(1) $, $(2)$, and $(3) $. Note also that the equivalence
$(1)\Leftrightarrow(4)$ is simply obtained by Lemma \ref{lem-GFD}
above and by the definition of the Gorenstein flat dimension.
Then, it remains to prove the equivalences
 $(1)\Leftrightarrow(2)\Leftrightarrow (3)$.\\
 $(1)\Rightarrow(2)$. We proceed by induction on $n$. We may assume
that $\Gfd_R(M)=n$. The case $n=0$ holds from Lemma
\ref{lem-close}. Then, suppose that $n\geq 1$. So, there exists,
by the definition of the Gorenstein flat dimension, a short exact
sequence of left $R$-modules:
$$ 0\rightarrow K\rightarrow G\rightarrow M \rightarrow 0,$$ where
$G$ is Gorenstein flat, such that $\Gfd_R(K)=n-1$. Then, for an
injective right $R$-module $I$, $\Tor_i^R( I,G)=0$ for all $i>0$
(by Lemma \ref{lem-close}), and  $\Tor_i^R(I,K )=0$ for all
$i>n-1$ (by induction). Then, we use the long exact sequence, $$
 \Tor_{i+1}^R(I,G )  \rightarrow  \Tor_{i+1}^R(I,M ) \rightarrow
 \Tor_{i}^R(I,K ),$$
to conclude that $\Tor_{i+1}^R(I,M)=0$ for all $i>n-1$, as
desired.\\
$(2)\Rightarrow(3)$. Easy by induction on $\id_R(E)$.\\
$(3)\Rightarrow(1)$.  Since $\Gfd_R(M)$ is finite and  by Lemma
\ref{lem-GFD}, we may pick, for some positive integer $m>n$, an
exact sequence of left $R$-modules:
$$0\rightarrow G_m\rightarrow \cdots \rightarrow G_0\rightarrow
M\rightarrow 0,$$ where $G_0,...,G_{m-1}$ are flat and $G_m$ is
Gorenstein flat. Let $K_n=\Ker(G_{n-1}\rightarrow G_{n-2})$. Our
aim is to prove that $K_n$ is Gorenstein flat. We decompose the
sequence
$$0\rightarrow G_m\rightarrow  \cdots \rightarrow G_n\rightarrow
K_n\rightarrow 0
$$ into short exact sequences:
$$ 0\rightarrow H_{i+1}\rightarrow      G_i\rightarrow
H_i\rightarrow 0$$ for $i=n,...,m-1$, where $H_n=K_n$ and
$H_{m}=G_{m}$. Consider the short exact sequence of left
$R$-modules:
$$0 \rightarrow H_{m}(=G_{m})\rightarrow   G_{m-1}\rightarrow
H_{m-1}\rightarrow 0.$$ We claim that $H_{m-1}$ is Gorenstein
flat. By the exact sequence:
$$0\rightarrow H_{m-1}\rightarrow  \cdots \rightarrow
G_0\rightarrow M\rightarrow 0, $$ we have, for every injective
right $R$-module $I$ and every $i>0$:
$$\Tor^R_{i}(I,H_{m-1} ) \cong \Tor^R_{(m-1)+i}(I,M )=0.$$
Therefore, by   Theorem \ref{thm-close} $(1)\Leftrightarrow(3)$,
 $ H_{m-1}$ is Gorenstein flat.\\
Finally, we repeat successively this last argument to conclude
that $ H_{m-2},...,H_{n}=K_n$ are Gorenstein flat, This completes
the proof.\cqfd\bigskip

Next generalizes \cite[Proposition 3.13]{HH} which is proved over
right coherent rings using the connection that exists between
Gorenstein flat dimension and Gorenstein injective dimension
\cite[Proposition 3.11]{HH}.

\begin{prop}\label{prop-sum-Gfd} If $R$ is  a left GF-closed ring, then
for a family of left $R$-modules $(M_i)_{i\in I}$, we have:
$$\Gfd_R(\oplus_{i\in I} M_i)=\sup\{\Gfd_R(M_i)\, |\, i\in I\}.$$
\end{prop}
\proof Using Theorem \ref{thm-GFD}, \cite[Proposition 3.2]{HH},
and the fact that over a left GF-closed ring $R$ the class
$\GF(R)$ is closed under direct summands (Corollary
\ref{cor-GF-closed}), the proof is analogous to the one of
\cite[Proposition 2.19]{HH}.\cqfd\bigskip

We end with the following generalization of \cite[Theorem
3.15]{HH}. Namely, we extend the standards (in)equalities for the
flat dimension; compare to \cite[Corollary 2, p. 135]{Bou}.

\begin{thm}\label{thm-short-Gfd}  If $R$ is  a left GF-closed ring, then
for  a short exact sequence of left $R$-modules  $0\rightarrow
A\rightarrow B\rightarrow C \rightarrow 0$, we have:
\begin{enumerate}
       \item If any two of the modules $A$, $B$, or $C$ have
finite Gorenstein flat dimension, then so has the third.
  \item  $\Gfd_R(A)\leq \sup\{\Gfd_R(B),\Gfd_R(C)-1\}$ with equality if
    $\Gfd_R(B)\not= \Gfd_R(C)$.
    \item $\Gfd_R(B)\leq \sup\{\Gfd_R(A),\Gfd_R(C)\}$ with equality if
    $\Gfd_R(C)\not= \Gfd_R(A)+1$.
    \item  $\Gfd_R(C)\leq \sup\{\Gfd_R(B),\Gfd_R(A)+1\}$ with equality if
    $\Gfd_R(B)\not= \Gfd_R(A)$.
\end{enumerate}
\end{thm}
\proof First note that, using the statement $(1)$ and Theorem
\ref{thm-GFD}, the statements $(2), (3),$ and $(4)$ are proved
similarly to \cite[Corollary 2, p. 135]{Bou}.\\
We prove $(1)$. Let $ \cdots \rightarrow A_1 \rightarrow
A_0\rightarrow A\rightarrow 0 $ and $ \cdots \rightarrow C_1
\rightarrow C_0\rightarrow C\rightarrow 0 $ be projective
resolutions of, respectively, $A$ and $C$. Then, by the Horseshoe
lemma \cite[Lemma 6.20]{Rot}, we get the following commutative
diagram, for any positive integer $n$:
$$\begin{array}{ccccccccc}
    &   & 0 &   & 0 &   & 0 &  &   \\
     &   & \downarrow &   & \downarrow  &   & \downarrow &   &  \\
   0 & \rightarrow & H_n & \rightarrow & K_n& \rightarrow&L_n  &\rightarrow & 0 \\
       &   & \downarrow &   & \downarrow  &   & \downarrow &   &  \\
     &   & \vdots &   & \vdots  &   & \vdots &   &  \\
   &   & \downarrow &   & \downarrow  &   & \downarrow &   &  \\
   0 & \rightarrow & A_0 & \rightarrow & A_0\oplus C_0 & \rightarrow&C_0 & \rightarrow &0 \\
     &   & \downarrow &   & \downarrow  &   & \downarrow &   &  \\
  0 & \rightarrow & A & \rightarrow &B& \rightarrow&C &\rightarrow & 0 \\
     &   & \downarrow &   & \downarrow  &   & \downarrow &   &  \\
     &   & 0 &   & 0 &   & 0 &  &
  \end{array}
$$
There are three cases:
\begin{description}
    \item[Case 1:] Assume that $\Gfd_R(A)\leq n$ and $\Gfd_R(C)\leq n$ for some positive
integer $n$. Then, by Theorem \ref{thm-GFD}, $H_n$ and $L_n$ are
Gorenstein flat, hence so is $K_n$ (since $R$ is left GF-closed).
This means that $\Gfd_R(B)\leq n$, as desired.
    \item[Case 2:]  Assume that $\Gfd_R(B)\leq n$ and $\Gfd_R(C)\leq n$ for some positive
integer $n$. So $K_n$ and $L_n$ are Gorenstein flat, hence so is
$H_n$ (by Theorem \ref{thm-close}). This means that $\Gfd_R(A)\leq
n$, as desired.
    \item[Case 3:] Assume that $\Gfd_R(A)\leq n$ and $\Gfd_R(B)\leq n$ for some positive
integer $n$. Then, $H_n$ and $K_n$ are Gorenstein flat. Assembling
the sequences  $ 0\rightarrow L_{n} \rightarrow C_{n-1}
\rightarrow
 \cdots \rightarrow C_0\rightarrow C\rightarrow 0 $ and $ 0\rightarrow H_{n}
\rightarrow K_{n} \rightarrow L_{n} \rightarrow 0 $, we get the
following exact sequence:  $$0\rightarrow H_{n} \rightarrow K_{n}
\rightarrow C_{n-1} \rightarrow  \cdots \rightarrow C_0\rightarrow
C\rightarrow 0.$$ Therefore, $\Gfd_R(C)\leq n+1$, as desired.\cqfd
\end{description}

\end{section}\bigskip\bigskip

\begin{section}{Gorenstein flat dimension in direct products of
rings}  Our aim, in this section, is to show, via a study of the
Gorenstein flat dimension of left modules over a direct product of
rings, how to construct a left GF-closed ring  that is neither
right coherent nor of finite weak dimension.\bigskip

For the convenience of the reader, we first recall some properties
concerning the structure of modules and homomorphisms over direct
products of rings (for more details please see \cite[Section
2.6]{BK}).\bigskip

Let $R=\displaystyle\prod_{i=1}^n R_i$ be a direct product of
rings. If $M_i$ is a left (resp., right) $R_i$-module for
$i=1,...,n$, then $M=M_1\oplus\cdots\oplus M_n$ is a left (resp.,
right) $R$-module. Conversely, if $M$ is  a left (resp., right)
$R$-module, then it is of the form $M=M_1\oplus\cdots\oplus M_n$,
where  $M_i$ is  a left (resp., right) $R_i$-module for
$i=1,...,n$ \cite[Subsection 2.6.6]{BK}. Also, the homomorphisms
of $R$-modules are determined by their actions on the $R_i$-module
components.  This is summarized in the following result:

\begin{thm}[\cite{BK}, Theorem 2.6.8]\label{thm-Hom-structure-product} Let
$R=\displaystyle\prod_{i=1}^n R_i$ be a direct product of rings
and let $M=M_1 \oplus \cdots\oplus M_n$ and
$N=N_1\oplus\cdots\oplus N_n$ be decompositions of  left (resp.,
right) $R$-modules into left (resp., right)  $R_i$-modules. Then,
the following hold:
\begin{enumerate}
    \item There is a natural isomorphism of abelian groups:
 $$\begin{array}{ccl}
   \Hom_R(M,N)&\stackrel{\cong}\longrightarrow   &  \Hom_{R_1}(M_1,N_1)\oplus
   \cdots
    \oplus\Hom_{R_n}(M_n,N_n)  \\
    \alpha & \longmapsto & \alpha_1\oplus \cdots
    \oplus \alpha_n \\
 \end{array}
$$
where the homomorphism $\alpha_1\oplus \cdots
    \oplus \alpha_n$ is defined by: $$(\alpha_1\oplus \cdots
    \oplus \alpha_n)(m_1,...,m_n)=(\alpha_1m_1, ...,
   \alpha_n m_n).$$
    \item The homomorphism $\alpha$ is injective (resp.,
    surjective) if  and only if  each $\alpha_i$ is injective (resp.,
    surjective).
\end{enumerate}
\end{thm}

Similarly, the tensor product of modules over a finite direct
product of rings is determined  as follows:

\begin{prop}\label{pro-tensor-structure-product}  Let
$R=\displaystyle\prod_{i=1}^n R_i$ be a direct product of rings
and let $M=M_1 \oplus \cdots\oplus M_n$   (resp.,
$N=N_1\oplus\cdots\oplus N_n$) be a decomposition of a right
(resp., left) $R$-module into right (resp., left)  $R_i$-modules.
Let $x=(x_1,..., x_n)$ and $y=(y_1,..., y_n)$ be elements of $M$
and $N$, respectively. Then, there is a natural isomorphism of
abelian groups:
 $$\begin{array}{ccl}
   M\otimes_R N &\stackrel{\cong}\longrightarrow   &( M_1\otimes_{R_1} N_1)\oplus \cdots
    \oplus ( M_n\otimes_{R_n}  N_n )  \\
   x\otimes_R y  & \longmapsto &   (x_1 \otimes_{R_1} y_1 ,...,x_n \otimes_{R_n} y_n).
 \end{array}
$$
\end{prop}

This may be used to give a proof of the following result.

\begin{prop}\label{prop-flat-product}
Let $R$ and $N$ be as in Proposition
\ref{pro-tensor-structure-product}. Then,
$$\fd_R(N)=\sup\{\fd_{R_i}(N_i)\, |\, i=0,...,n \}.$$
\end{prop}

Now, we state the main result in this section.

\begin{thm}\label{thm-Gflat-product}
Let $R$ and $N$ be as in Proposition
\ref{pro-tensor-structure-product}. Then,
$$\Gfd_R(N)=\sup\{\Gfd_{R_i}(N _i)\, |\, i=0,...,n \}.$$
\end{thm}
\proof First, as a consequence of Theorem
\ref{thm-Hom-structure-product}, Propositions
\ref{pro-tensor-structure-product} and \ref{prop-flat-product},
and the fact that a right $R$-module  $E=E_1 \oplus \cdots\oplus
E_n$ is injective if  and only if  each $E_i$ is an injective
right  $R_i$-module (by \cite[Exercise 2.6.5]{BK}), we get the
equivalence:  $N $ is a Gorenstein  flat left $R$-module if and
only if each $N _i$ is a Gorenstein flat left $R_i$-module.\\
Now, we prove the desired equality.  First, we prove the inequality
$\Gfd_R(N )\leq \sup\{\Gfd_{R_i}(N _i)\, |\, i=0,...,n \}.$ For
that, we may assume that $\sup\{\Gfd_{R_i}(N _i)\, |\, i=0,...,n
\} = m$ for some positive integer  $m$. Then, there exists, for
$i=1,...,n$, an exact sequence of left $R_i$-modules:
$$     0\longrightarrow F_{m,i} \stackrel{  \alpha_{m,i}}
\longrightarrow\cdots\longrightarrow   F_{0,i}\stackrel{
\alpha_{0,i}}\longrightarrow  N _{i}\longrightarrow 0, $$ where $
F_{j,i}$ is  Gorenstein  flat   for $j=0,..., m$. Then, we get an
exact sequence of left $R$-modules:
$$     0\longrightarrow F_{m} \stackrel{  \alpha_{m}}
\longrightarrow\cdots\longrightarrow   F_{0}\stackrel{
\alpha_{0}}\longrightarrow  N \longrightarrow 0, $$ where $F_{j}=
F_{j,1}\oplus \cdots \oplus F_{j,n}$ for $j=0,..., m$. Since every
$ F_{j,i}$ is a  Gorenstein  flat  left $R_i$-module, every
$F_{j}$ is a Gorenstein  flat left $R$-module (by
the reason above). This implies the first inequality.\\
Now, we prove the converse inequality. For that, we may assume
that $\Gfd_R(N )= m$ for some positive integer $m$. By Theorem
\ref{thm-Hom-structure-product}, we have an exact sequence of left
$R$-modules:
 $$    0 \longrightarrow \oplus_i F_{m,i} \stackrel{ \oplus_i\alpha_{m,i}}
\longrightarrow\cdots \longrightarrow \oplus_i F_{0,i}\stackrel{
\oplus_i\alpha_{0,i}}  \longrightarrow N = \oplus_i N
_{i}\longrightarrow 0,$$ where  $\oplus_i F_{j,i}$ is a Gorenstein
flat left $R$-module for $j=1,..., m$. So, by the reason above,
each $F_{j,i}$ is a Gorenstein  flat left $R_i$-module for
$i=1,..., n$ and $j=1,..., m$ . Then, for $i=1,..., n$, we have an
exact sequence of  left $R_i$-modules:
$$     0\longrightarrow F_{m,i} \stackrel{  \alpha_{m,i}}
\longrightarrow\cdots\longrightarrow   F_{0,i}\stackrel{
\alpha_{0,i}}\longrightarrow  N _{i}\longrightarrow 0, $$ where $
F_{j,i}$ is   Gorenstein  flat for $j=0,..., m$. Thus, for
$i=1,..., n$, $\Gfd_{R_i}(N _i)\leq m$. Therefore, $
\sup\{\Gfd_{R_i}(N _i)\, |\, i=0,...,n \}\leq m=\Gfd_R(N ).$ This
completes the proof.\cqfd\bigskip

Now, we are in position to give the desired example. For that,  we
need the following result.

\begin{prop}\label{pro-product-GF-closed} A direct product of rings
$\displaystyle\prod_{i=1}^n R_i$ is a left GF-closed ring if  and
only if each $R_i$ is left GF-closed.
\end{prop}
\proof Simply combine Theorem \ref{thm-Hom-structure-product} with
Theorem \ref{thm-Gflat-product}.\cqfd

Finally, the results above give as a method to construct examples
of left  GF-closed rings which are  neither right coherent nor of
finite weak dimension, as follows:

\begin{exmp}\label{exm-3.2} \textnormal{From the note after Proposition
\ref{pro-example-GF-closed}, we may have a left GF-closed ring $R$
with infinite weak dimension and a left GF-closed ring  $S$ which
is not right coherent. Then, for the direct product of rings
$\Gamma=R\times S$ we have:
\begin{enumerate}
 \item $\Gamma$ is left GF-closed (by Proposition \ref{pro-product-GF-closed}).
     \item $\wdim(\Gamma)=\sup\{ \wdim(R), \wdim(S)\}=\infty$
     (for the left equality use Proposition \ref{prop-flat-product}).
    \item $\Gamma$ is not right coherent (since  using Theorem
    \ref{thm-Hom-structure-product} we may prove that a product of
    rings $R_1\times R_2$ is right coherent if and only if
   each $R_i$ is right coherent\,\footnote{\; This equivalence is a
   particular case, in commutative setting, of
   \cite[Theorem 2.13]{DKM}.}.)
\end{enumerate}}
\end{exmp}



\noindent {\bf Acknowledgement.} The author is grateful to the
referee for the valuable comments and suggestions.

\end{section}



\bigskip\bigskip

\end{document}